\documentclass[letterpaper,11pt]{article}

 \title{A note on linear B-spline copulas}
 \author{Arturo Erdely}
 \date{\small{Facultad de Estudios Superiores Acatl\'an \\
              Universidad Nacional Aut\'onoma de M\'exico \\
							\texttt{arturo.erdely@comunidad.unam.mx}\\}}
 
 \usepackage{amsthm}
 \usepackage{amssymb}
 \usepackage{amsmath}
 \usepackage{enumerate}
 \usepackage{graphicx}
 
 \newcommand{\indic}{\textbf{\textsf{\large{1}}}}
 
\begin{document}
 
\maketitle

\begin{abstract}
  \noindent In this brief note we prove that linear B-spline copulas is not a new family of copulas since they are equivalent to checkerboard copulas, and discuss in particular how they are used to extend empirical subcopulas to copulas.
\end{abstract}


\section{Introduction}

In the context of probability and statistics, copula functions are used to represent the functional relationship between the joint distribution function of a random vector and its corresponding univariate marginal distribution functions. In the 2-dimensional case, for example, if $(X,Y)$ is a random vector and $F_{XY}$ its joint distribution function, the outstanding theorem by Sklar (1959) may be used in this particular case to prove that there exists a 2-copula function $C$ such that $F_{XY}(x,y)=C(F_X(x),F_Y(y))$ where $F_X$ and $F_Y$ are the marginal univariate distribution functions of $X$ and $Y,$ respectively. A \textit{2-copula} is a function $C:I^2\rightarrow I,$ where $I:=[0,1],$ such that $C(u,0)=0=C(0,v),$ $C(u,1)=u$ and $C(1,v)=v,$ and that is \textit{2-increasing}: the $C$-volume $V_C(B)\geq 0$ for any rectangle $B=[u_1,u_2]\times[v_1,v_2]\subseteq I^2$ where $V_C(B)=C(u_2,v_2)-C(u_2,v_1)-C(u_1,v_2)+C(u_1,v_1).$ Copulas are always uniformly continuous.

\bigskip

\noindent In fact, \textit{Sklar's Theorem} proves that there is a unique \textit{subcopula} that links a joint distribution function to its univariate margins, and thereafter it is proved that such unique subcopula can be extended to a generally nonunique copula, but the copula is unique if the univariate margins are continuous, see for example Schweizer and Sklar (1983) or Nelsen (1999). The definition of \textit{2-subcopula} is almost the same of a copula except for the fact that its domain is a subset $S_1\times S_2$ of the unit square $I^2$ where $\{0,1\}\subseteq S_i\subseteq I,$ and therefore a 2-copula is a 2-subcopula whose domain is $I^2.$

\bigskip

\noindent We will prove that \textit{linear B-spline copulas} as in Shen \textit{et al.} (2008) is not a new family of copulas since they are the same as the ones obtained by bilinear extension of certain type of subcopulas, defined as \textit{checkerboard copulas} by Li \textit{et al.} (1997) or Li et al. (1998), and discuss its application to extend empirical subcopulas to copulas.

\section{Checkerboard 2-copulas}

\noindent Let $Q_1^{(m_1)}$ and $Q_2^{(m_2)}$ be two subsets of $I$ such that $Q_k^{(m_k)}=\{q_{ki}:i=0,\ldots,m_k\}$ where $m_k\in\{1,2,\ldots\}$ and $0=q_{k0}<q_{k1}<\cdots<q_{k,m_k-1}<q_{km_k}=1$ for $k\in\{1,2\}.$ If $C:I^2\rightarrow I$ is a 2-copula, then the restriction of $C$ to the mesh $D:=Q_1^{(m_1)}\times Q_2^{(m_2)}$ is a subcopula which will be denoted as $S.$ That is, for any $(u,v)\in D$ we have $S(u,v)=C(u,v).$

\bigskip

\noindent Define subsets $Q_{ij}\subseteq I^2$ such that $Q_{ij}:=\langle q_{1,i-1},q_{1i}]\times\langle q_{2,j-1},q_{2j}]$ for $i\in\{1,\ldots,m_1\}$ and $j\in\{1,\ldots,m_2\},$ with the convention that $\langle$ is for left-open if $i$ or $j$ is equal to $1$ and left-closed otherwise. Then $\mathbb{Q}:=\{Q_{ij}:i=1,\ldots,m_1\text{ and }j=1,\ldots,m_2\}$ is a partition of $I^2.$ For any $(u,v)\in Q_{ij}$ define the functions
\begin{equation}
\lambda_{1i}(u) := \,\frac{u - q_{1,i-1}}{q_{1i}-q_{1,i-1}} \qquad\qquad \lambda_{2j}(v) := \,\frac{v - q_{2,j-1}}{q_{2j}-q_{2,j-1}} 
\end{equation}
The function $C_{\mathbb{Q}}:I^2\rightarrow I$ defined for each $(u,v)$ in some $Q_{ij}$ by
\begin{eqnarray}\label{bilinearCopula}
   C_{\mathbb{Q}}(u,v) &:=& [1-\lambda_{1i}(u)][1-\lambda_{2j}(v)]S(q_{1,i-1},q_{2,j-1}) + [1-\lambda_{1i}(u)]\lambda_{2j}(v)S(q_{1,i-1},q_{2j}) \nonumber \\
			                 &{ }& + \lambda_{1i}(u)[1-\lambda_{2j}(v)]S(q_{1i},q_{2,j-1}) + \lambda_{1i}(u)\lambda_{2j}(v)S(q_{1i},q_{2j})
\end{eqnarray}
is a copula that extends sucopula $S$ to $I^2$ as proved, for example, by Lemma 2.3.5 in Nelsen (1999), that may be considered as an approximation of $C,$ and that in particular agrees with $S$ (and therefore with $C$) on $D.$ It is straightforward to obtain the copula density associated to $C_{\mathbb{Q}}$ which for each $(u,v)$ in some $Q_{ij}$ is given by
\begin{equation}\label{bilinearDensity}
  c_{\mathbb{Q}}(u,v) =\, \frac{\partial^2}{\partial u\partial v}C_{\mathbb{Q}}(u,v) =\, \frac{V_S(Q_{ij})}{V_{\Pi}(Q_{ij})}
\end{equation}
where $V_S(Q_{ij})$ is the $S-$volume of $Q_{ij}:$
\begin{equation}\label{SVolume}
   V_S(Q_{ij}) = S(q_{1i},q_{2j}) - S(q_{1,i-1},q_{2j}) - S(q_{1i},q_{2,j-1}) + S(q_{1,i-1},q_{2,j-1})
\end{equation}
and $V_{\Pi}(Q_{ij})$ the $\Pi-$volume of $Q_{ij}$ where $\Pi(u,v)=uv,$ that is:
\begin{equation}\label{PiVolume}
  V_{\Pi}(Q_{ij}) = q_{1i}q_{2j} - q_{1,i-1}q_{2j} - q_{1i}q_{2,j-1} + q_{1,i-1}q_{2,j-1} = (q_{1i} - q_{1,i-1})(q_{2j} - q_{2,j-1})
\end{equation}
It should be noticed that since $S$ is the restriction of $C$ to $D$ then $V_S(Q_{ij})=V_C(Q_{ij}).$

\bigskip

\noindent Now consider the particular case where $m_1=m_2=m\in\{1,2,\ldots\}$ and $Q^{(m)}=\{i/m: i=0,1,\ldots,m\}.$ If a subcopula $S$ is the restriction of a copula $C$ to the mesh $D=Q^{(m)}\times Q^{(m)}$ then the extension of subcopula $S$ to a copula $C_{\mathbb{Q}}$ as in (\ref{bilinearCopula}) is called a \textbf{checkerboard copula} approximation of $C,$ which in fact approximates $C$ uniformly, see Li \textit{et al.} (1997) or Li et al. (1998). For this particular case we obtain the following expression for the checkerboard copula density, where $\mathbf{1}_{Q_{ij}}$ is the indicator function of $Q_{ij}=\langle\frac{i-1}{m},\frac{i}{m}]\times\langle\frac{j-1}{m},\frac{j}{m}]:$
\begin{equation}\label{checkerboardDensity}
  c_{\mathbb{Q}}(u,v) = m^2\sum_{i\,=\,1}^m\sum_{j\,=\,1}^m V_C(Q_{ij})\mathbf{1}_{Q_{ij}}(u,v)
\end{equation}

\bigskip

\noindent Shen \textit{et al.}(2008) proposed a linear B-spline approximation $C_{LB}$ of a copula $C$ by the following:
\begin{equation}
  C_{LB}(u,v)=\sum_{i\,=\,0}^m\sum_{j\,=\,0}^m C\bigg(\frac{i}{m},\frac{j}{m}\bigg)N_i(u)N_j(v)
\end{equation}
where for $i=1,\ldots,m-1$
      $$N_i(u)=\left\{ \begin{array}{ll}
	                            mu - (i - 1), & \frac{i-1}{m}\leq u < \frac{i}{m}\,, \\
									            { }           & { } \\
									            (i + 1) - mu, & \frac{i}{m}\leq u < \frac{i+1}{m}\,, \\
									            { }           & { } \\
									            0,            & \text{otherwise,}
	                     \end{array}\right.$$

      $$N_0(u)=\left\{ \begin{array}{ll}
			                        1 - mu, & 0\leq u < \frac{1}{m}\,, \\
									            0,      & \text{otherwise,}
			                 \end{array}\right.$$
			$$N_m(u)=\left\{ \begin{array}{ll}
			                        mu - (m - 1), & \frac{m-1}{m} \leq u \leq 1, \\
									            0,      & \text{otherwise,}
			                 \end{array}\right.$$
and proved that $C_{LB}$ is a copula that uniformly approximates $C.$ It is straightforward to notice that the copula density for $C_{LB}$ is the same as (\ref{checkerboardDensity}):
\begin{equation}
  c_{LB}(u,v) = \frac{\partial^2}{\partial u\partial v}C_{LB}(u,v) = m^2\sum_{i\,=\,1}^m\sum_{j\,=\,1}^m V_C(Q_{ij})\mathbf{1}_{Q_{ij}}(u,v) = c_{\mathbb{Q}}(u,v)
\end{equation}
and therefore the linear B-spline copula approximation is the same as the checkerboard copula approximation for any given copula. This proves that the claim in Shen \textit{et al.}(2008) that linear B-spline copulas is a new family of copulas is false.

\section{Empirical 2-subcopula and extensions}

Let $\{(x_1,y_1),\ldots,(x_n,y_n)\}$ denote an observed bivariate sample of size $n$ from a continuous distribution. In Nelsen (1999) it is defined as \textit{empirical copula} the function $C_n$ given by
\begin{equation}\label{empiricalCopula}
  C_n\bigg(\,\frac{i}{n}\,,\,\frac{j}{n}\,\bigg)\,:=\,\frac{1}{n}\sum_{k\,=\,1}^n\indic\{\text{rank}(x_k)\leq i\,,\,\text{rank}(y_k)\leq j\}\,,\qquad i,j\in\{0,\ldots,n\}
\end{equation}
Strictly speaking, what is usually known as empirical copula is not a copula, it is a \textbf{subcopula} with domain the mesh $D=\{i/n:i=0,1,\ldots,n\}^2,$ which of course can be extended to a checkerboard copula as explained in the previous section. So what Shen \textit{et al.}(2008) define as \textit{empirical linear B-spline copula} is in fact a checkerboard copula extension of the empirical (sub)copula given in (\ref{empiricalCopula}), and such authors redefine ``empirical copula'' as a discontinuous extension of (\ref{empiricalCopula}) by a function $\hat{C}_E:I^2\rightarrow I$ given by
$$ \hat{C}_E(u,v) = \left\{ \begin{array}{ll}
                                  C_n\big(\,\frac{i-1}{n}\,,\,\frac{j-1}{n}\,\big), & \frac{i-1}{n}\leq u < \frac{i}{n}\,,\,\frac{j-1}{n} \leq v < \frac{j}{n}\,, \\
																	1,   & u = v = 1,
	                          \end{array} \right.$$
but being discontinuous $\hat{C}_E$ clearly is not a copula, in contrast with the checkerboard approximation which extends (\ref{empiricalCopula}) to a copula. Therefore, the comparison of nonparametric estimation of copulas performed by Shen \textit{et al.}(2008) between the empirical linear B-spline copula and their definition of empirical copula is, in fact, a comparison of the checkerboard copula extension of (\ref{empiricalCopula}) versus a discontinuous noncopula extension of (\ref{empiricalCopula}).

\section{Final remarks}

This note was submitted for publication to \textit{Computational Statistics \& Data Analysis} and received the following response from the reviewer:
\begin{quote}
 The manuscript correctly states that linear B-spline copulas as in Shen et al. (2008) are the same as the ones obtained by bilinear extension of certain type of subcopulas, defined as checkerboard copulas by Li et al. (1997). It was a pity that the 2008 paper did not explicitly state this link, however, nowadays, there have been various papers that have discussed checkerboard copulas, clarified possible links with the literature, studied possible applications to empirical copulas. To cite a few of them, please see:\par\smallskip
*Carley, H. and Taylor, M. (2003). A new proof of Sklar's theorem. In Cuadras, C. M., Fortiana, J., and Rodr\'iguez Lallena, J., editors, Distributions with given marginals and Statistical Modelling, pages 29-34. Kluwer, Dordrecht.\par
*Genest, C., Neslehova, J. G., and Remillard, B. (2014). On the empirical multilinear copula process for count data. Bernoulli, 20(3):1344-1371.\par
*Ghosh, S. and Henderson, S. G. (2009). Patchwork distributions. In Alexopoulos, C., Goldsman, D., and Wilson, J., editors, Advancing the frontiers of simulation:
a Festschrift in honor of George Samuel Fishman, pages 65-86. Springer.\par
*Johnson, N. L. and Kotz, S. (2007). Cloning of distributions. Statistics, 41(2):145-152.\par\smallskip
I do not think the present short note adds more insights into the problem and, hence, I do not think it is suitable for publication.
\end{quote}

\noindent As anyone can check in detail, all of the above citations have no reference to Shen \textit{et al.}(2008) nor a mention about the equivalence between checkerboard copulas and linear B-spline copulas. Moreover, a recently published authoritative book on copula theory by Durante and Sempi (2016) which includes a section on checkerboard copulas does not mention its equivalence to linear B-spline copulas, not even in the historical remarks about the subject, which by the way it does include the above reviewer's citations. Finally, in the recently published and also authoritative book on copula modeling by Joe (2015) bivariate linear B-spline copulas and citation of Shen \textit{et al.}(2008) is included as a non-parametric estimation method of low-dimensional copulas, but no mention at all about its equivalence with checkerboard copulas. Therefore, it is a pity that a prestigious journal such as \textit{Computational Statistics \& Data Analysis} refuses to publish a note which their reviewers accept is correct about a published article that included a false claim about a new class of copulas.

\section*{References}

\noindent Carley, H. and Taylor, M. (2003) A new proof of Sklar's theorem. In Cuadras, C. M., Fortiana, J., and Rodr\'iguez Lallena, J., editors, \textit{Distributions with given marginals and Statistical Modelling}, pages 29--34. Kluwer (Dordrecht).
\medskip

\noindent Durante, F. and Sempi, C. (2016) \textit{Principles of Copula Theory.} CRC Press (Boca Raton).
\medskip

\noindent Genest, C., Ne\v{s}lehov\'a, J. G., and R\'emillard, B. (2014) On the empirical multilinear copula process for count data. \textit{Bernoulli} \textbf{20} (3), 1344--1371.
\medskip

\noindent Ghosh, S. and Henderson, S. G. (2009) Patchwork distributions. In Alexopoulos, C., Goldsman, D., and Wilson, J., editors, \textit{Advancing the frontiers of simulation:
a Festschrift in honor of George Samuel Fishman}, pages 65--86. Springer.
\medskip

\noindent Joe, H. (2015) \textit{Dependence Modeling with Copulas.} CRC Press (Boca Raton).
\medskip

\noindent Johnson, N. L. and Kotz, S. (2007) Cloning of distributions. \textit{Statistics} \textbf{41} (2), 145--152.
\medskip

\noindent Li, X., Mikusi\'nski, P., Sherwood, H., and Taylor, M.D. (1997) On approximation of copulas. In Bene\v{s}, V. and \v{S}t\v{e}p\'an, J., editors, \textit{Distributions with given marginals and moment problems,} pp. 107--116. Kluwer (Dordrecht).
\medskip

\noindent Li, X., Mikusi\'nski, P., and Taylor, M.D. (1998) Strong approximation of copulas. \textit{J. Math. Anal. Appl.} \textbf{225}, 608--623.
\medskip

\noindent Nelsen, R.B. (1999) \textit{An Introduction to Copulas.} Springer (New York).
\medskip

\noindent Schweizer, B. and Sklar, A. (1983) \textit{Probabilistic Metric Spaces.} North-Holland (New York).
\medskip

\noindent Shen, X., Zhu, Y., and Song, L. (2008) Linear B-spline copulas with applications to nonparametric estimation of copulas. \textit{Computational Statistics \& Data Analysis} \textbf{52}, 3806--3819. doi: 10.1016/j.csda.2008.01.002
\medskip

\noindent Sklar, A. (1959) Fonctions de r\'epartition \`a $n$ dimensions et leurs marges. \textit{Publ. Inst. Statist. Univ. Paris,} \textbf{8}, 229--231.

\end{document}